\documentstyle[12pt]{article}

\def\be{\begin{equation}}
\def\ee{\end{equation}}
\def\bea{\begin{eqnarray}}
\def\eea{\end{eqnarray}}

\def\1{\'{\i}}                           
\def\R{{\bf R}}

\def\kk{K}
\def\pp{P}
\def\hh{H}
\def\cuno{C_1}
\def\cdos{C_2}
\def\dd{D}
\def\zzt{\tau}

\begin{document}
\begin{center}
\vspace*{1.0cm}

{\LARGE{\bf A new quantum so(2,2) algebra}}

\vskip 1.25cm

{\large {\bf Francisco J. Herranz\footnote{
Contribution to the
Proceedings of the `Quantum Theory and Symmetries'
(Goslar, 18-22 July 1999) (World Scientific, 2000),
edited by H.-D. Doebner, V.K. Dobrev, J.-D. Hennig and W.
Luecke.} }}

\vskip 0.5 cm

Departamento de F\1sica\\
Escuela Polit\'ecnica Superior\\ 
Universidad de Burgos\\
E--09006 Burgos, Spain

\end{center}

\vspace{0.75 cm}

\begin{abstract}
By starting from the non-standard quantum deformation of the $sl(2,\R)$
algebra, a new quantum deformation for the real Lie algebra $so(2,2)$
is constructed by imposing the former to be a Hopf subalgebra of the
latter. The quantum $so(2,2)$ algebra so obtained is realized as a
quantum conformal algebra of the $(1+1)$  Minkowskian
spacetime. This Hopf algebra is  shown to be the symmetry
algebra of a  time discretization of the $(1+1)$  wave
equation and its   contraction gives rise to a new  $(2+1)$ quantum  
  Poincar\'e algebra.
\end{abstract}

\vspace{0.75 cm}


\section{Introduction}

The non-standard  quantum deformation of $sl(2,\R)\simeq so(2,1)$
\cite{Ohn}, here denoted $U_z(sl(2,\R))$, has been the starting point in
the obtention of  non-standard  quantum algebras in higher dimensions.
In particular, by taking two copies of $U_z(sl(2,\R))$ and applying
the same procedure as in the standard (Drinfel'd--Jimbo) case
\cite{Ita}, a quantum
 $so(2,2)$ algebra has been obtained in \cite{beyond}, while the
corresponding deformation for $so(3,2)$ has been found in \cite{vulpi}.
These  quantum algebras have been realized as deformations of
conformal algebras for the Minkowskian spacetime. Furthermore, by
following either a contraction approach \cite{beyond} or a deformation
embedding method \cite{null}, non-standard quantum deformations for
other Lie algebras have been deduced; amongst them  it is remarkable
the appearance of a non-standard quantum Poincar\'e algebra,
which can be considered  as  a  conformal quantum algebra for
the Carroll spacetime, or alternatively  
as a null-plane quantum Poincar\'e algebra  \cite{null}.  All these
results are summarized in the following diagram where the vertical
arrows indicate the corresponding contractions leading  to Poincar\'e
algebras:
$$
\begin{array}{ccccc}
U_z(sl(2,\R))&\longrightarrow & U_z(sl(2,\R))\oplus
U_{-z}(sl(2,\R))\simeq U_z(so(2,2))&\longrightarrow &
U_z(so(3,2))\\[5pt]
 \Big\downarrow\
\varepsilon
\to 0&&
\qquad\qquad\qquad\qquad\qquad\qquad\quad 
\Big\downarrow\ \varepsilon
\to 0& & \Big\downarrow\ \varepsilon \to 0\\[5pt]
U_z(iso(1,1))&\longrightarrow&
\quad\mbox{\footnotesize{Null-plane  Poincar\'e algebra}}\qquad
U_z(iso(2,1))&\longrightarrow & U_z(iso(3,1))\cr
 \end{array} 
$$

The aim of this contribution is to provide, starting  again from
$U_z(sl(2,\R))$, a new way in the obtention of non-standard quantum
algebras. The first step is to construct a new non-standard quantum
$so(2,2)$ algebra which could be the cornerstone of further constructions
in higher dimensions.  The essential idea is to require that
$U_z(sl(2,\R))$ remains as a Hopf subalgebra so that this approach can
be seen as a kind of {\em complete} deformation embedding method.
Next, a contraction limit  gives rise to a new $(2+1)$ 
quantum Poincar\'e algebra which  contains a $(1+1)$ 
quantum Poincar\'e Hopf subalgebra:
$$
U_z(sl(2,\R))\subset U_z(so(2,2))
\quad\stackrel{\varepsilon \to 0}{\longrightarrow}\quad
 U_z(iso(1,1))\subset U_z(iso(2,1))
$$

It is interesting to stress that such new  quantum $so(2,2)$
algebra   is   the symmetry algebra of a time
discretization of the wave equation. Thus we recall in the next section
the basic facts of the Lie algebra $so(2,2)$ in a conformal basis as well
as its relationship with the $(1+1)$   wave equation. The Hopf
algebra structure deforming  $so(2,2)$, its role as a discrete symmetry
algebra and its contraction to Poincar\'e are  presented in the section
3.


\section{Lie algebra so(2,2)}

Let us consider the real Lie algebra $so(2,2)$ generated by
  $H$   (time translations), $P$  (space
translations), $K$ (boosts), $D$ (dilations) and $C_1$, $C_2$ (special
conformal transformations). In this basis 
$so(2,2)$ is the Lie algebra of the group of  conformal
transformations of the  $(1+1)$  Minkowskian spacetime.
 The Lie brackets of $so(2,2)$ read
\be
\begin{array}{lll}
[K,H]=P&\qquad [K,P]=H&\qquad [H,P]=0\cr
[D,H]=H&\qquad [D,C_1]=-C_1&\qquad [H,C_1]=-2D\cr
[D,P]=P&\qquad [D,C_2]=-C_2&\qquad [P,C_2]=2D\cr
[K,C_1]=C_2&\qquad [K,C_2]=C_1&\qquad [C_1,C_2]=0\cr
[H,C_2]=2K&\qquad [P,C_1]=-2K&\qquad [K,D]=0 .
 \end{array} 
 \label{ba}
\ee
Three subalgebras of $so(2,2)$ are relevant for our purposes:\\
$\bullet$ $\{H,P,K\}$ which span the $(1+1)$  Poincar\'e
algebra (first row in (\ref{ba})). 

\noindent
$\bullet$ $\{D,H,C_1\}$ which give rise to  $so(2,1)\simeq sl(2,\R)$
(second row in (\ref{ba})).

\noindent
$\bullet$ $\{D,P,C_2\}$ which also generate  $so(2,1)\simeq sl(2,\R)$
(third row in (\ref{ba})).

A vector field  representation of $so(2,2)$ in terms of the space and
time coordinates $(x,t)$ is given by
\be
\begin{array}{l}
 H=\partial_t\qquad P=\partial_x\qquad
K =- t \partial_x -x\partial_t\qquad
D=-x\partial_x - t \partial_t\cr
C_1=(x^2+t^2) \partial_t + 2 x t   \partial_x\qquad
C_2=-(x^2+t^2) \partial_x - 2 x t   \partial_t .
 \end{array} 
 \label{bb}
\ee
The Casimir of the above Poincar\'e subalgebra is $E=P^2-H^2$. The
action of $E$ on a function $\Phi(x,t)$ through the 
representation (\ref{bb}) (choosing for $E$ the value zero) leads to
the $(1+1)$  wave equation:
\be
E\Phi(x,t)=0\qquad
\Longrightarrow\qquad\left(\frac{\partial^2}{\partial
x^2}-\frac{\partial^2}{\partial t^2}
\right)\Phi(x,t)=0 .
\label{bc}
\ee
We shall say that an operator ${\cal O}$ is a symmetry of the
equation $E\Phi(x,t)=0$ if ${\cal O}$ transforms solutions into
solutions, that is, $E{\cal O}=\Lambda E$ where $\Lambda$ is another
operator. The Lie algebra $so(2,2)$ is the symmetry algebra of the
wave equation: $E$ commutes with $\{H,P,K\}$ and  in the
realization (\ref{bb}) the remaining  generators  verify
\be
 [E,D]=-2E\qquad
[E,C_1]= 4 t E\qquad [E,C_2]= - 4 x E .
 \label{bd}
\ee


\section{Non-standard quantum so(2,2) algebra}

We choose the  $sl(2,\R)$ subalgebra of $so(2,2)$ spanned by 
$\{D,H,C_1\}$. Then we write in terms of these generators the
non-standard  quantum deformation of $sl(2,\R)$ in the form introduced
in \cite{BH} and denote $\tau$  the deformation parameter. This  means
that the classical $r$-matrix we are considering for $so(2,2)$ is 
$r=-\zzt\dd\wedge \hh$ (which is a solution of the classical
Yang--Baxter equation). Now we look for a quantum 
$so(2,2)$ algebra that keeps the quantum $sl(2,\R)$ algebra as a Hopf
subalgebra: $U_\tau(sl(2,\R))\subset U_\tau(so(2,2))$.  The resulting
coproduct and commutation rules for $U_\tau(so(2,2))$ are given by:
\be
\begin{array}{l}
\Delta(\hh)=1\otimes \hh + \hh \otimes 1\qquad
\Delta(\pp)= 1\otimes \pp + \pp \otimes e^{\zzt\hh}\cr
\Delta(\dd)= 1\otimes \dd + \dd \otimes e^{-\zzt\hh}
\qquad\Delta(\cuno)= 1\otimes \cuno + \cuno \otimes
e^{-\zzt\hh}\cr
\Delta(\kk)=1\otimes \kk + \kk\otimes 1 - 
\zzt \dd \otimes e^{-\zzt\hh}\pp\cr
\Delta(\cdos)=1\otimes \cdos + \cdos
\otimes e^{-\zzt\hh}  + 2 \zzt\dd \otimes e^{-\zzt\hh}\kk
-\zzt^2\dd(\dd+1)\otimes e^{-2\zzt\hh}\pp
\end{array} 
 \label{cb}
\ee
\be
\begin{array}{l}
[\kk,\hh]=e^{-\zzt\hh}\pp\qquad 
[\kk,\pp]= ( e^{\zzt\hh}-1)/{\zzt} 
\qquad [\hh,\pp]=0\cr
 [\dd,\hh]=({1- e^{-\zzt\hh}})/{\zzt}
\qquad [\dd,\cuno]=-\cuno +   \zzt \dd^2 \qquad
[\hh,\cuno]=-2\dd\cr
 [\dd,\pp]= \pp\qquad [\dd,\cdos]=-\cdos \qquad
[\pp,\cdos]=2\dd\cr
 [\kk,\cuno]=\cdos  \qquad
[\kk,\cdos]=\cuno- \zzt \dd^2\qquad
[\cuno,\cdos]=-\zzt(\dd\cdos+\cdos\dd)\cr
 [\hh,\cdos]=e^{-\zzt\hh}\kk + \kk
e^{-\zzt\hh}\quad\ [\pp,\cuno]=-2\kk - 
\zzt(\dd\pp+\pp\dd)\quad\ [\kk,\dd]=0 .
\end{array} 
 \label{cc}
\ee
It can be checked that the universal quantum $R$-matrix for 
$U_\tau(sl(2,\R))$ \cite{BH} also holds for $U_\tau(so(2,2))$. In our
basis this element reads
\be
{\cal R}=\exp\left\{\tau H\otimes D\right\}
\exp\left\{-\tau D\otimes H\right\} .
\label{cd}
\ee

The relationship between $U_\tau(so(2,2))$ and a discretization of the 
wave equation can be established by means of the following 
differential-difference realization which under
the limit $\tau\to 0$ gives the classical realization (\ref{bb}):
\bea
&&\hh=\partial_t \qquad \pp=\partial_x\nonumber\\[2pt]
&&\kk= - x\left( \frac{ e^{\zzt\partial_t}-1
}{\zzt}\right) - t  e^{-\zzt\partial_t}  \partial_x\qquad
\dd =- x \partial_x  
- t \left( \frac{1-e^{-\zzt\partial_t}}{\zzt}\right)   \cr 
&&\cuno=(x^2 +t^2 e^{- \zzt\partial_t})\left( \frac{
e^{\zzt\partial_t}-1}{\zzt}\right)
+ 2 x t \partial_x +\zzt  x \partial_x  
+ \zzt  x^2 \partial^2_{x}\cr
&&\cdos=-(x^2+t^2e^{-2\zzt\partial_t})\partial_x
 - 2 x t \left( \frac{1- e^{-\zzt\partial_t} }{\zzt} \right)
+\zzt t  e^{-2 \zzt\partial_t}\partial_x .
\label{da}
\eea
The generators  $\{\hh,\pp,\kk\}$ close  a deformed Poincar\'e
subalgebra (although not a Hopf subalgebra) whose Casimir is now 
$E_\zzt=\pp^2 -\left( \frac{e^{\zzt\hh}-1}{\zzt}\right)^2$.
If we introduce the realization (\ref{da}) then we find a time
discretization of the wave equation on a uniform lattice with $x$ as a
continuous variable:
\be
E_\zzt\Phi(x,t)=0\qquad
\Longrightarrow\qquad \left\{
\frac{\partial^2}{\partial x^2}-
\left(\frac{e^{\zzt\partial_t}-1}{\zzt}\right)^2 
\right\}\Phi(x,t)=0 .
\label{dc}
\ee
Therefore the deformation parameter $\tau$ appearing within the 
discrete derivative in (\ref{dc}) can be identified with the time
lattice constant. Furthermore the generators (\ref{da}) are symmetry
operators of (\ref{dc}) since they fullfil
\bea
&& [E_\zzt,X]=0\quad \mbox{for}\quad X\in\{\hh,\pp,\kk\}
\qquad
[E_\zzt,\dd]=-2E_\zzt \cr
&&[E_\zzt,\cuno]= 4(t +\zzt + \zzt x \partial_x) E_\zzt
\qquad [E_\zzt,\cdos]= - 4 x E_\zzt .
 \label{dd}
\eea
Hence  we conclude that $U_\zzt(so(2,2))$ is the symmetry algebra of
the discrete wave equation (\ref{dc}). In this respect we recall
that the symmetries of a discretization of the  wave equation in both
coordinates $(x,t)$ on a uniform lattice were computed in
\cite{Luismi},  showing that they are  difference   operators which
preserve the Lie algebra $so(2,2)$ as in the continuous case. Therefore
some kind of connection between the  results of \cite{Luismi} and our 
quantum
$so(2,2)$ algebra should exist  as it was already
established for discrete Shr\"odinger equations and quantum algebras
\cite{schrod}.

To end with, we work out the contraction from $U_\zzt(so(2,2))$ to a
new quantum Poincar\'e algebra: $U_\zzt(so(2,2))\to U_\zzt(iso(2,1))$.
We apply to the Hopf algebra $U_\zzt(so(2,2))$ the In\"on\"u--Wigner
contraction defined by the map
\be
H\to \varepsilon H \qquad P\to  P \qquad K\to \varepsilon K \qquad 
C_1\to \varepsilon C_1 \qquad C_2\to  C_2 \qquad D\to D
\label{de}
\ee
together with a transformation of the deformation parameter:
$\tau\to \tau/\varepsilon$. The limit $\varepsilon\to 0$ leads to the
coproduct and commutators of $U_\zzt(iso(2,1))$:
\be
\begin{array}{l}
\Delta(\hh)=1\otimes \hh + \hh \otimes 1\qquad
\Delta(\pp)= 1\otimes \pp + \pp \otimes e^{\zzt\hh}\cr
\Delta(\dd)= 1\otimes \dd + \dd \otimes e^{-\zzt\hh}\qquad
\Delta(\cuno)= 1\otimes \cuno + \cuno \otimes
e^{-\zzt\hh} \cr
\Delta(\kk)=1\otimes \kk + \kk\otimes 1\qquad
\Delta(\cdos)=1\otimes \cdos + \cdos
\otimes e^{-\zzt\hh}  + 2 \zzt\dd \otimes e^{-\zzt\hh}\kk
\end{array} 
 \label{df}
\ee
\be
\begin{array}{lll}
[K,H]=0&\qquad [K,P]=(e^{\zzt\hh}-1)/{\zzt} &\qquad
[H,P]=0\cr [D,H]=({1- e^{-\zzt\hh}})/{\zzt}&\qquad
[D,C_1]=-C_1&\qquad [H,C_1]=0\cr [D,P]=P&\qquad [D,C_2]=-C_2&\qquad
[P,C_2]=2D\cr [K,C_1]=0&\qquad [K,C_2]=C_1&\qquad [C_1,C_2]=0\cr
[H,C_2]=2 e^{-\zzt\hh}\kk&\qquad [P,C_1]=-2K&\qquad [K,D]=0 .
 \end{array} 
 \label{dg}
\ee
The universal quantum $R$-matrix for $U_\zzt(iso(2,1))$ is also 
(\ref{cd}) so that it is formally preserved  under
contraction. Note also that the generators $\{D,H,C_1\}$ give rise to a
$(1+1)$ quantum Poincar\'e subalgebra such that:
 $U_\zzt(iso(1,1))\subset  U_\zzt(iso(2,1))$. 
 
Finally we remark that if  we would have chosen the 
$sl(2,\R)$ subalgebra spanned by $\{D,P,C_2\}$ instead of the 
one generated by $\{D,H,C_1\}$, then we would have obtained a quantum
$so(2,2)$ algebra with $P$ as primitive generator (instead of $H$). This
second choice would lead to a space discretization of the wave equation.  
Both quantum $so(2,2)$ algebras would be
algebraically equivalent by the interchanges $H\leftrightarrow P$ and 
$C_1\leftrightarrow C_2$, however their contraction would lead to
inequivalent  quantum Poincar\'e algebras. A complete  analysis
of all these possibilities will be presented elsewhere.


\section*{Acknowledgment}
 This work was partially supported  by  Junta de Castilla y
Le\'on, Spain  (Project   CO2/399).


\end{document}